# Potential symmetries and conservation laws for generalized quasilinear hyperbolic equations

M. Nadjafikhah*    R. Bakhshandeh-Chamazkoti†    F. Ahangari‡

*School of Mathematics, Iran University of Science and Technology, Narmak, Tehran 1684613114, Iran.*

**Abstract**

Based on Lie group method, potential symmetry and invariant solutions for generalized quasilinear hyperbolic equations are studied. To obtain the invariant solutions in explicit form, we focus on the physically interesting situations which admit potential symmetries. Then by using the partial Lagrangian approach, we find conservation laws for this equation in three physically interesting cases.

**MSC 2010:** 70S10, 35L65, 70H33
**Keywords:** Conservation laws, Generalized quasilinear hyperbolic equations, Invariant solution, Potential symmetry.

## 1 Introduction

Lie groups involving potential symmetries is exposed in view of applications to physics. Also, the concept of a conservation law and the relationship between symmetries and conservation laws, which are mathematical formulations of the familiar physical laws of conservation of energy, conservation of momentum and so on, plays an important role in the analysis of basic properties of the solutions, [1, 2, 3, 4, 5].

In 1988, Bluman et al. [6] suggested a method for obtaining a new class of symmetries for a partial differential equation **R**, which is written in a conservative form. They analyzed the Lie symmetries of the system **S** that is obtained by introducing a potential as a further unknown function. In this approach, we obtain classes of symmetries by computing Lie groups of point transformations, whose infinitesimals act on a different space than the space of independent, dependent variables, and their derivatives of the system **S**. These symmetries are neither point symmetries nor Lie-Bäcklund symmetries. The existence of potential symmetries leads to the construction of corresponding invariant solutions. Any group of Lie transformations for the system **S** induces a potential symmetry for **R**, when at least one of the generators depends explicitly on the potential. The potential symmetries are useful in looking for the solutions using a reduction method, [2, 7, 8].

A systematic way for the determination of conservation laws associated with variational symmetries for systems of Euler-Lagrange equations is indeed the famous Noether theorem, [9, 10]. This theorem requires a Lagrangian. There are approaches that do not require a Lagrangian or even assume the existence of a Lagrangian for differential equations, [11, 12]. A new idea is presented by Kara and Mahomed in [13] on how one can construct conservation laws of differential equations via operators that are not necessarily symmetry generators of the underlying system. These partial Noether operators which are associated with what are called

---

*m_nadjafikhah@iust.ac.ir
†r_bakhshandeh@iust.ac.ir
‡fa_ahangari@iust.ac.ir

*partial Lagrangians*, aid via an explicit formula in the construction of conservation laws of the underlying system which need not be derivable from a variational principle. These systems are called *partial Euler-Lagrange equations* with respect to partial Lagrangians. This method provides a systematic manner of obtaining conservation laws for systems which have partial Lagrangian formulations, [14, 15].

The main purpose of the present paper is to find potential symmetries and conservation laws of the generalized quasilinear hyperbolic equations by using the nonlocal conservation theorem method and the partial Lagrangian approach. In addition, invariant solution is obtained by utilizing the relationship between conservation laws and Lie-point symmetries of the equation

$$f(x)u_{tt} = [g(x,u)u_x + h(x,u)]_x, \tag{1}$$

where $f \in \mathcal{C}^1(\mathbb{R})$ and $g, h \in \mathcal{C}^1(\mathbb{R}^2)$ are arbitrary nonzero smooth functions. In case of $f(x) = 1$ and $g_x = h_x = 0$ the equation (1) is changed to nonlinear telegraph equation for which, the point symmetries and nonlocal potential symmetries are obtained by George W. Bluman et al. in [16].

The organization of this paper is as follows:

In Section 2, we present fundamental definitions and theorems about potential symmetry and nonlocal conservation laws. In Section 3, we discuss the potential symmetries of the generalized quasilinear hyperbolic equation (1) and we calculate these symmetries for three interesting cases. In section 4, We obtain the invariant solutions associated to the potential symmetries discussed in the previous section. Section 5 is devoted to the study of the conservation laws of the equation (1). Concluding remarks are given in the final Section.

## 2 Preliminaries

Consider the $\mathbf{R}\{x, t; u\}$ is a partial differential equation (PDE) of $n$ independent variables $x = (x^1, \ldots, x^n)$ with components $x^i$ and one dependent variable $u$ with components $u^\alpha$ by

$$R(x, u, u_{(1)}, \ldots, u_{(r)}) = 0. \tag{2}$$

where $u_{(1)}, \ldots, u_{(r)}$ are the collections of all first, second, ..., $k$th order partial derivatives, i.e., $u_i^\alpha = D_i(u^\alpha), u_{ij}^\alpha = D_j D_i(u^\alpha), \ldots$ respectively, with the total derivative operator with respect to $x^i$ given by

$$D_i := \frac{\partial}{\partial x^i} + u_i^\alpha \frac{\partial}{\partial u^\alpha} + u_{ij}^\alpha \frac{\partial}{\partial u_j^\alpha} + \ldots, \quad i = 1, \ldots, n. \tag{3}$$

Let $\mathcal{A}$ be the space of differential functions. The Lie-Bäcklund or generalized operator is

$$X := \xi^i \frac{\partial}{\partial x^i} + \eta^\alpha \frac{\partial}{\partial u^\alpha} + \sum_{s \geq 1} \zeta_{i_1 \ldots i_s}^\alpha \frac{\partial}{\partial u_{i_1 \ldots i_s}^\alpha}, \quad \xi^i, \eta^\alpha \in \mathcal{A} \tag{4}$$

where $\zeta_{i_1 \ldots i_s}^\alpha$ are defined by

$$\zeta_i^\alpha = D_i(\eta^\alpha) - u_j^\alpha D_i(\xi^j) \tag{5}$$

$$\zeta_{i_1 \ldots i_s}^\alpha = D_{i_s}(\eta_{i_1 \ldots i_s}^\alpha) - u_{j i_1 \ldots i_s}^\alpha D_{i_s}(\xi^j) \tag{6}$$

The Noether operators associated with a Lie-Bäcklund operator $X$ are

$$N^i = \xi^i + W^\alpha \frac{\partial}{\partial u_i^\alpha} + \sum_{s \geq 1} D_{i_1} \ldots D_{i_s} \frac{\partial}{\partial u_{i_1 \ldots i_s}^\alpha}, \tag{7}$$

in which $W^\alpha$ is the Lie characteristic function by definition

$$W^\alpha = \eta^\alpha - \xi^j u_j^\alpha, \tag{8}$$

and the Euler-Lagrange operator $\frac{\delta}{\delta u_i^\alpha}$ is

$$\frac{\delta}{\delta u_i^\alpha} = \frac{\partial}{\partial u_i^\alpha} + \sum_{s \geq 1}(-1)^s D_{j_1} \ldots D_{j_s} \frac{\partial}{\partial u_{j_1 \ldots j_s}^\alpha}, \quad i = 1, 2, \ldots, n; \ \alpha = 1, 2, \ldots, r. \tag{9}$$

The $n$-tuple vector

$$T = (T^1(x, u, u_{(1)}, \ldots, u_{(r-1)}), \ldots, T^n(x, u, u_{(1)}, \ldots, u_{(r-1)})),$$

is a *conserved vector* of (2) if $D_i(T_i)\big|_{(2)} = 0$; This equation defines a *local conservation law* of system (2). George W. Bluman et al. [6] have introduced the concept of potential symmetry for any differential equation which can be written as a conservative form

$$D_t F(x, t, u, u_x, u_t) - D_x G(x, t, u, u_x, u_t) = 0, \tag{10}$$

where $D_x$ and $D_t$ are the total derivative operators. Introducing an auxiliary potential variable $v = v(x, t)$, it is possible to form the potential system, $S = 0$,

$$v_x = F, \qquad v_t = G. \tag{11}$$

The potential system $\mathbf{S}\{x, t; u, v\}$ (11) has essentially the same solution set as that of the PDE $\mathbf{R}\{x, t; u\}$ (10), [8].

**Theorem 2.1** *equation (10) admits symmetry potentials only if equation (11) assumes one of the $v_x = H(x, t, u)u_t + K_1(x, t, u)u_x + K_2(x, t, u)$ forms where $H \neq 0$; otherwise $v_x = K(x, t, u, u_x)$ and in this case $\tau = \tau(t)$.*

Let the system (2) is written in the form

$$E_\alpha = E_\alpha^0 + E_\alpha^1 = 0, \quad \alpha = 1 \ldots, m. \tag{12}$$

If there exists a function $L = L(x, u, u_{(1)}, \ldots, u_{(r)}) \in \mathcal{A}$, $l \leq k$ and nonzero functions $f_\alpha^\beta \in \mathcal{A}$ such that (12) can be written as

$$\frac{\delta L}{\delta u^\alpha} = f_\alpha^\beta E_\beta^1, \quad \alpha, \beta = 1, \ldots, m. \tag{13}$$

provided $E_\beta^1 = 0$ for some $\beta$, $L$ is called a *partial Lagrangian* of (12) otherwise it is the standard Lagrangian. We term differential equations of the form (13) *partial Euler–Lagrange equations*.

A Lie-Bäcklund or generalized operator $X$ of the form (4) is called a *partial Noether operator* corresponding to a partial Lagrangian $L \in \mathcal{A}$ if and only if there exists a vector $B = (B^1, \ldots, B^n)$, $B^i \in \mathcal{A}$ such that

$$X(L) + LD_i(\xi^i) = W^\alpha \frac{\delta L}{\delta u^\alpha} + D_i(B^i), \tag{14}$$

where $W = (W^1, \ldots, W^m)$, $W^\alpha \in \mathcal{A}$, is the characteristic of $X$, [4].

**Theorem 2.2** *(Partial Noether Theorem) A Lie–Bäcklund operator $X$ of the form (4) is a partial Noether operator of a partial Lagrangian $L$ corresponding to the partial Euler–Lagrange system (13) if and only if the characteristic $W = (W, \ldots, W)$, $W \in \mathcal{A}$, of $X$ is also the characteristic of the conservation law $D_i T^i = 0$, where*

$$T^i = B^i - N^i(L), \quad i = 1, \ldots, n. \tag{15}$$

*of the partial Euler–Lagrange equations (13).*

## 3 Potential symmetries

In order to find the potential symmetries of equation (1), we set it in the following conservative form:

$$D_t\big(f(x)u_t\big) - D_x\big(g(x,u)u_x + h(x,u)\big) = 0, \tag{16}$$

By considering a potential $v(x,t)$ as an auxiliary unknown function, the following system can be associated with (16)

$$v_x = f(x)u_t, \qquad v_t = g(x,u)u_x + h(x,u), \tag{17}$$

The system, which characterizes the generators is obtained from

$$X^{(1)}\big(v_x - fu_t\big)\Big|_{(3.17)} = 0, \qquad X^{(1)}\big(v_t - gu_x - h\big)\Big|_{(3.17)} = 0, \tag{18}$$

which must hold identically. Here, $X^{(1)}$ is the operator

$$X^{(1)} = \xi\frac{\partial}{\partial x} + \tau\frac{\partial}{\partial t} + \varphi\frac{\partial}{\partial u} + \eta\frac{\partial}{\partial v} + \varphi_1^{(1)}\frac{\partial}{\partial u_x} + \varphi_2^{(1)}\frac{\partial}{\partial u_t} + \eta_1^{(1)}\frac{\partial}{\partial v_x} + \eta_2^{(1)}\frac{\partial}{\partial v_t}, \tag{19}$$

where

$$\begin{aligned}
\varphi_1^{(1)} &= \varphi_x + (\varphi_u - \xi_x)u_x - \tau_x u_t - \tau_u u_x u_t - \xi_u u_x^2 + \varphi_v v_x - \tau_v u_t v_x - \xi_v u_x v_x, \\
\varphi_2^{(1)} &= \varphi_t + (\varphi_u - \tau_t)u_t - \tau_u u_t^2 - \xi_t u_x - \xi_u u_t u_x + \varphi_v v_t - \tau_v u_t v_t - \xi_v u_x v_t, \\
\eta_1^{(1)} &= \eta_x + (\eta_v - \xi_x)v_x - \xi_v v_x^2 + \eta_u u_x - \tau_x v_t - \tau_u u_x v_t - \tau_v v_x v_t - \xi_u u_x v_x, \\
\eta_2^{(1)} &= \eta_t + (\eta_v - \tau_t)v_t - \tau_v v_t^2 + \eta_u u_t - \tau_u u_t v_t - \xi_t v_x - \xi_u u_t v_x - \xi_v v_t v_x.
\end{aligned}$$

Therefore the (18) becomes

$$\begin{aligned}
\eta_x &+ (\eta_v - \xi_x)v_x - \xi_v v_x^2 + (\eta_u + f\xi_t)u_x - (\tau_x + f\varphi_v)v_t + (f\xi_v - \tau_u)u_x v_t - \tau_v v_x v_t \\
&- \xi_u u_x v_x - f(\varphi_t - \tau_u u_t^2 - \xi_u u_t u_x - \tau_v u_t v_t) - [f_x\xi + f(\varphi_u - \tau_t)]u_t = 0, \tag{20}\\
\eta_t &+ (\eta_v - \tau_t)v_t - \tau_v v_t^2 + (\eta_u + g\tau_x)u_t - \tau_u u_t v_t - (\xi_t + \varphi_v)v_x + (\tau_v - \xi_u)u_t v_x - \xi_v v_t v_x \\
&- g[\varphi_x - \tau_u u_x u_t - \xi_u u_x^2 - \xi_v u_x v_x] + [g_x\xi - g_u\varphi - g(\varphi_u - \xi_x)]u_x + h_u\varphi + h_x\xi = 0.
\end{aligned}$$

On substituting $v_x$ by $f(x)u_t$, and $v_t$ by $g(x,u)u_x + h(x,u)$ in equations (20) and (21), we get the determining equations

$$\tau_u = f\xi_v, \tag{21}$$
$$\xi_u = g\tau_v, \tag{22}$$
$$\eta_u - gf\varphi_v + h\tau_u = 0, \tag{23}$$
$$\eta_x - f\varphi_t - fh\varphi_v - h\tau_x = 0, \tag{24}$$
$$f(\eta_v + \tau_t - \varphi_u - \xi_x) - f_x\xi = 0, \tag{25}$$
$$h_x\xi - \eta_t + h\tau_t - h\eta_v + h_u\varphi + g\varphi_x + h^2\tau_v = 0, \tag{26}$$
$$g(\xi_x - \tau_t - \varphi_u + \eta_v) - g_u\varphi - g_x\xi - 2\xi_u h = 0. \tag{27}$$

In solving the above system of equations (21)–(27), we confine our attention to the relevant physical forms of $f(x)$, $g(x,u)$ and arbitrary smooth function $h(x,u)$.

The point symmetry group for (17) is completely determined by the generators $\xi$, $\tau$, $\varphi$ and $\eta$. Those point symmetries which satisfy $\xi_v^2 + \tau_v^2 + \varphi_v^2 = 0$, correspond to point symmetries of (1). Instead, the potential symmetries of (1) are obtained, if $\xi_v^2 + \tau_v^2 + \varphi_v^2 > 0$. This group maps a solution of (17) into another solution of (17), thus induces a mapping of a solution of (1) into another solution of (1).

In the following, we will mention some cases which admits potential symmetries.

## 3.1 $f(x) = a$, $g(x,u) = u$ and $h(x,u) = \exp(x)$ where $a$ is constant.

In this case, the infinitesimal symmetries are given by the following:

$$\xi = c_1 x + c_2, \qquad \tau = c_1 t + c_3,$$
$$\varphi = c_1 u + c_4 v + c_6 t + c_5, \quad \eta = c_1 v + c_4 a\left(\frac{u^2}{2} + e^x\right) + c_6 a x + c_7.$$

Then, we obtain point symmetries with the following generators:

$$X_1 := x\frac{\partial}{\partial x} + t\frac{\partial}{\partial t} + u\frac{\partial}{\partial u} + v\frac{\partial}{\partial v}, \qquad X_2 := \frac{\partial}{\partial x}, \qquad X_3 := \frac{\partial}{\partial t},$$

$$X_4 := v\frac{\partial}{\partial u} + a\left(\frac{u^2}{2} + e^x\right)\frac{\partial}{\partial v}, \qquad X_5 := \frac{\partial}{\partial u}, \qquad X_6 := ax\frac{\partial}{\partial v}, \qquad X_7 := \frac{\partial}{\partial v}.$$

Clearly, $X_4$ is a potential symmetry for (1), since we have: $\xi_v^2 + \tau_v^2 + \varphi_v^2 = 1 > 0$.

## 3.2 $f(x) = x$, $g(x,u) = xu$ and $h(x,u) = u^2$.

In this case, the infinitesimal symmetries are

$$\xi = c_1 x + c_2, \qquad \tau = c_1 t + c_3,$$
$$\varphi = c_1 u + c_4(2vt + x^2 u) + 2c_5 t + 2c_6 v + c_8,$$
$$\eta = c_1 v + c_4 x^2(v + u^2 t) + c_5 x^2 + c_6 x^2 u^2 + c_7.$$

Therefore, we have point symmetries with the following generators

$$X_1 := x\frac{\partial}{\partial x} + t\frac{\partial}{\partial t} + u\frac{\partial}{\partial u} + v\frac{\partial}{\partial v}, \qquad X_2 := \frac{\partial}{\partial x},$$

$$X_3 := \frac{\partial}{\partial t}, \qquad X_4 := (2vt + x^2 u)\frac{\partial}{\partial u} + x^2(v + u^2 t)\frac{\partial}{\partial v},$$

$$X_5 := 2t\frac{\partial}{\partial u} + x^2\frac{\partial}{\partial v}, \qquad X_6 := 2v\frac{\partial}{\partial u} + x^2 u^2\frac{\partial}{\partial v},$$

$$X_7 := \frac{\partial}{\partial v}, \qquad X_8 := \frac{\partial}{\partial u}.$$

It is clear that $X_4$ and $X_6$ are only potential symmetries for (1). Since we have $\xi_v^2 + \tau_v^2 + \varphi_v^2 = 2t > 0$ and $\xi_v^2 + \tau_v^2 + \varphi_v^2 = 2 > 0$, respectively.

## 3.3 $f(x) = a$, $g(x,u) = u$ and $h(x,u) = u$, where $a$ is an arbitrary constant.

In this case, we have the infinitesimal symmetries

$$\xi := c_1 v + c_3 x + c_2,$$
$$\tau := c_1 a(x + u) + c_3 t + c_4,$$
$$\varphi := -c_1 v + c_3 u + c_4(x + u) + c_5(t + 1),$$
$$\eta := (c_3 + c_4)v + c_5 a x + c_6.$$

with following generators

$$X_1 := v\frac{\partial}{\partial x} + a(x+u)\frac{\partial}{\partial t} - v\frac{\partial}{\partial u}, \qquad X_2 := \frac{\partial}{\partial x},$$

$$X_3 := x\frac{\partial}{\partial x} + t\frac{\partial}{\partial t} + u\frac{\partial}{\partial u} + v\frac{\partial}{\partial v}, \qquad X_4 := \frac{\partial}{\partial t} + (x+u)\frac{\partial}{\partial u} + v\frac{\partial}{\partial v},$$

$$X_5 := (t+1)\frac{\partial}{\partial u} + ax\frac{\partial}{\partial v}, \qquad X_6 := \frac{\partial}{\partial v}.$$

Obviously, $X_1$ is the only potential symmetry for (1). Since, we have: $\xi_v^2 + \tau_v^2 + \varphi_v^2 = 2 > 0$.

# 4 Invariant solutions

Given a point symmetry for (17), the invariant surface conditions are:

$$\xi(x,t,u,v)u_x + \tau(x,t,u,v)u_t - \varphi(x,t,u,v) = 0,$$
$$\xi(x,t,u,v)v_x + \tau(x,t,u,v)v_t - \eta(x,t,u,v) = 0.$$

The solutions of the associated characteristic system are defined as one-parameter families of characteristic curves [7].

## 4.1 For case of 3.1

For the potential symmetry $X_4$ from the invariant surface condition we infer that:

$$v = 0, \quad u(x,t) = \sqrt{2\exp(x)}i.$$

## 4.2 For case of 3.2

For the potential symmetry $X_4$, from the invariant surface condition we infer that:

$$u(x,t) = \frac{1}{2t^2 x^2}, \quad v = \frac{1}{4t^3}.$$

Taking into account the invariant surface condition of potential symmetry $X_6$ we obtain only trivial solution.

## 4.3 For case of 3.3

For the potential symmetry $X_1$ the invariant surface condition is as follows:

$$vu_x + a(x+u)u_t + v = 0$$
$$vv_x + a(x+u)v_t = 0$$

one of the solutions of the above system is

$$v = 0, \quad u = F(x)$$

by substituting $u = F(x)$ in equation (1) we obtain the following reduced ODE:

$$F'^2 + FF'' + F' = 0.$$

# 5 Conservation laws

Here, we further investigate conservation laws of the equation (1), in three possible cases for $f(x)$, $g(x,u)$ and $h(x,u)$ which were considered in Section 4.

The equation (1) has the partial Lagrangian form $L = \frac{1}{2}g(x,u)u_x^2 - \frac{1}{2}f(x)u_t^2$ and the Euler–Lagrange type equation is $\frac{\delta L}{\delta u} = \frac{1}{2}g_u u_x^2 + h_x + h_u u_x$ for which the partial Noether operators $X := \tau \partial_t + \xi \partial_x + \varphi \partial_u$ satisfy in (14), therefore we have

$$\begin{aligned}
\frac{1}{2}(\xi g_x + \varphi g_u)u_x^2 &- \frac{1}{2}\xi f' u_t^2 - f(x)u_t[\varphi_t + \varphi_u u_t - u_t(\tau_t + \tau_u u_t) - u_x(\xi_t + \xi_u u_t)] \\
&+ g(x,u)u_x[\varphi_x + \varphi_u u_x - u_t(\tau_x + \tau_u u_x) - u_x(\xi_x + \xi_u u_x)] \\
&+ (\tau_t + \tau_u u_t + \xi_x + \xi_u u_x)\left(\frac{1}{2}g(x,u)u_x^2 - \frac{1}{2}f(x)u_t^2\right) \\
&= (\varphi - \tau u_t - \xi u_x)\left(\frac{1}{2}g_u u_x^2 + h_x + h_u u_x\right) + B_t^1 + B_u^1 u_t + B_x^2 + B_u^2 u_x.
\end{aligned}$$
(28)

By solving (28), the partial Noether operator obtained is $X = \alpha(t,x)\frac{\partial}{\partial u}$ with gauge terms

$$B^1 = -f(x)\alpha_t(t,x)u + \beta(t,x), \quad B^2 = \alpha_x(t,x)\int g(x,u)du - \alpha(x,t)h(x,u) + \gamma(t,x), \quad (29)$$

where $\alpha(t,x)$ and $\beta(t,x)$ are arbitrary smooth functions with respect to $t, x$. By using (15) the conserved vector $T = (T^1, T^2)$ is obtained with the following conserved components

$$T^1 = -f(x)\alpha_t(t,x)u + \alpha(x,t)u_t + \beta(t,x), \quad (30)$$

$$T^2 = \alpha_x(t,x)\int g(x,u)du - \alpha(x,t)h(x,u) - \alpha(x,t)g(x,t)u_x + \gamma(t,x). \quad (31)$$

Now, we obtain the conserved vector $T = (T^1, T^2)$ in three subsections 3.1, 3.2 and 3.3 of the previous section. The results are presented as follows:

## 5.1 For case of 3.1

If $f(x) = a$, $g(x,u) = u$ and $h(x,u) = \exp(x)$ where $a$ is constant, so

$$T = \left(-a\alpha_t(t,x)u + \alpha(x,t)u_t + \beta(t,x),\ \alpha_x(t,x)\frac{u^2}{2} - \alpha(x,t)(\exp(x) + uu_x) + \gamma(t,x)\right). \quad (32)$$

## 5.2 For case of 3.2

If $f(x) = x$, $g(x,u) = xu$ and $h(x,u) = u^2$, then

$$T = \left(-x\alpha_t(t,x)u + \alpha(x,t)u_t + \beta(t,x),\ \alpha_x(t,x)\frac{xu^2}{2} - \alpha(x,t)(u^2 + xuu_x) + \gamma(t,x)\right). \quad (33)$$

## 5.3 For case of 3.3

If $f(x) = a$, $g(x,u) = u$ and $h(x,u) = u$, where $a$ is an arbitrary constant, we have

$$T = \left(-a\alpha_t(t,x)u + \alpha(x,t)u_t + \beta(t,x),\ \alpha_x(t,x)\frac{u^2}{2} - \alpha(x,t)(u + uu_x) + \gamma(t,x)\right). \quad (34)$$

# Conclusion

In this paper, we have given a complete analysis on finding potential symmetries for the generalized quasilinear hyperbolic equations. The infinitesimals, similarity variables, dependent variables, and reduction to quadrature or exact solutions of the mentioned generalized quasilinear hyperbolic equations for physically realizable forms of $f(x), g(x,u)$ and $h(x,u)$ are also obtained. Further conservation laws are constructed for this equation in three different cases. It will be interesting to study potential, nonclassical potential and nonclassical symmetries of generalized quasilinear hyperbolic equations to search for new exact solutions.